\documentclass[11pt]{amsart}
\usepackage{geometry}
\usepackage{array}
\usepackage{amssymb,amsthm,amsmath}
\usepackage{ulem}
\usepackage{amsmath}
\usepackage{moreverb}
\usepackage[all]{xy}
\usepackage{graphicx}

\usepackage{mathtools}
\usepackage{enumitem}
\usepackage[backend=bibtex,citestyle=alphabetic,bibstyle=alphabetic,maxbibnames=99]{biblatex}
\usepackage{hyperref}
\usepackage[compress]{cleveref}

\providecommand{\phantomsection}{}
\AtBeginDocument{\let\textlabel\label}
\makeatletter
\newcommand{\mylabel}[2]{\raisebox{.7\normalbaselineskip}{\phantomsection}#1%
  \def\@currentlabel{#1}\textlabel{#2}}
\newcommand{\mylabelrev}[2]{\raisebox{.7\normalbaselineskip}{\phantomsection}%
  \def\@currentlabel{#1}\textlabel{#2}}
\makeatother

\newcounter{eqsign}
\renewcommand{\theeqsign}{(\arabic{eqsign})}
\newcommand{\eqprintnew}[1]{\refstepcounter{eqsign} \overset{\theeqsign}{=} \mylabelrev{\theeqsign}{#1}}

{\theoremstyle{definition}  }
{}
{\theoremstyle{remark} }

{\theoremstyle{remark} }
\newtheorem{mainthm}{Theorem}

\newcommand{\CC}{\mathbb{C}} 
 
\newcommand{\NN}{\mathbb{N}} 
 
\newcommand{\schat}[1]{\mathfrak{S}_{#1}}

\newcommand{\disp}[0]{\displaystyle}
\DeclarePairedDelimiter{\bracepair}{\lbrace}{\rbrace}
\DeclarePairedDelimiter{\anglepair}{\langle}{\rangle}

\DeclarePairedDelimiter{\parenpair}{(}{)}
\DeclarePairedDelimiter{\vertpair}{\vert}{\vert}
\DeclarePairedDelimiter{\Vertpair}{\Vert}{\Vert}

\newcommand{\tra}{\mathop{\mathit{Trace}}}

\newcommand{\declare}{:=}

\DeclareMathOperator{\ranky}{Rank}

\makeatletter
\newcommand{\raisemath}[1]{\mathpalette{\raisem@th{#1}}}
\newcommand{\raisem@th}[3]{\raisebox{#1}{$#2#3$}}
\makeatother

\newcommand{\zd}[1]{\mathbb{Z}_{#1}}
\newcommand{\mcy}[1]{m \parenpair{#1}}
\newcommand{\spectr}[1]{\sigma \parenpair{#1}}
\newcommand{\snorm}[2]{\Vertpair{#1}_{#2}}
\newcommand{\sspace}[2]{\mathcal{N} \parenpair{#1 ; #2}}
\newcommand{\locrad}[1]{\delta (#1)}
\newcommand{\locproj}[1]{P(#1)}

\begin{document}
\title[Criterion for \texorpdfstring{$\zd{\MakeLowercase{d}}$}{Zd}--symmetry]{Criterion for \texorpdfstring{$\zd{\MakeLowercase{d}}$}{Zd}--symmetry of a Spectrum of a Compact Operator}
\author{Boris S. Mityagin}
\address{231 West 18th Avenue, The Ohio State University, Columbus, OH 43201}
\email{\url{borismit@math.ohio-state.edu}}

\begin{abstract}
If $A$ is a compact operator in a Banach space and some power $A^q$ is nuclear we give a criterion of $\zd{d}$ -- symmetry of its spectrum $\spectr{A}$ in terms of vanishing of the traces $\tra A^n$ for all $n$, $n \geq 0$, $n \neq 0 \mod d$, sufficiently large.
\end{abstract}
\maketitle
\normalem

In the case of matrices, or linear operators $T: X \to X$ in a finite-dimensional space, one can check (prove) that the following conditions are equivalent.
\begin{enumerate}[label = (\alph*)]
\item \label{enum:specsym} The spectrum of $T$ is symmetric, or $\zd{2}$-symmetric, i.e., $\lambda \in \spectr{T} \rightarrow - \lambda \in \spectr{T}$ and their algebraic multiplicities $\mcy{\lambda}, \mcy{-\lambda}$ are equal; 
\item \label{enum:matrtrace} $\tra T^p = 0$ for all odd $p \in \NN$.
\end{enumerate}
M. Zelikin \cite{Zelikin} observed and proved that this claim could be extended to $\schat{1}$, the trace-class operators in a Hilbert space.  We will show that such claims could be made 
\begin{enumerate}[label = (\roman*)]
\item in general Banach spaces; 
\item for $\zd{d}$ symmetry of a spectrum, $d \geq 2$.  
\end{enumerate}

Of course, we need to make sure that $\tra$ is well--defined if we write conditions like \ref{enum:matrtrace}.  Then, the formula for the trace $\tra A = \sum_{j} \lambda_j(A)$ should be properly explained if we use it.  We now recall a few notions and facts about nuclear operators (see more in \cite{Konigb}). 

An operator $A: X \to Y$ between two Banach spaces is called \textit{nuclear} if it has representation
\begin{equation} \label{eq:decomp}
Ax = \sum_{k = 1}^q a_k f_k(x) y_k, \quad q \leq \infty
\end{equation}
where 
\begin{equation} \label{eq:norms}
\begin{split}
a_k > 0, \quad a^* = \sum a_k < \infty, \quad \text{ and } \\ 
\Vertpair{f_k \vert X^{\prime}} \leq 1, \quad \Vertpair{y_k \vert Y} \leq 1, \quad \forall k.
\end{split}
\end{equation}
A linear space of nuclear operators $X \to Y$ is a Banach space $\sspace{X}{Y}$ with the norm
\begin{equation} \label{eq:sonenorm}
\snorm{A}{1} = \inf \bracepair{a^*: \eqref{eq:decomp}, \eqref{eq:norms}}.
\end{equation}
A linear functional $\tra$ is well-defined on $\sspace{X}{X}$ (for any Banach Space $X$) by 
\begin{equation} \label{eq:tradef}
\tra A = \sum_{k = 1}^q a_k f_k(y_k).
\end{equation}
Of course,
\begin{equation}
\Vertpair{\tra A} \leq \snorm{A}{1},
\end{equation}
and $\Vertpair{\tra} = 1$.  

A. Grothendieck \cite{Gro} showed that for operators \eqref{eq:decomp} $[X = Y]$,
\begin{align}
\text{ if } &\sum_{k = 1}^q a_k^{2/3} < \infty \label{eq:twothirds}\\
\text { then } & \sum \vertpair*{\lambda_j(A)} < \infty
\end{align}
where points of the spectrum $\spectr{A}$ are enumerated with their multiplicity, and 
\begin{equation} \label{eq:algsptreqiv}
\tra A = \sum \lambda_{j}(A).
\end{equation}
The presentation \eqref{eq:tradef} with \labelcref{eq:decomp,eq:norms,eq:sonenorm} gives a factorization
\begin{equation} \label{eq:factor}
A = JF, \quad X \xrightarrow{F} \ell_2(\NN) \xrightarrow{J} X,
\end{equation}
where 
\begin{align}
Fx &= \sum_{1}^{\infty} a_k^{1/2} f_k(x) e_k, \quad \text{ and }  \label{eq:fdef}\\
J \xi & = \sum_{1}^{\infty} a_k^{1/2} \xi_k y_k, \label{eq:jdef} 
\end{align}
with
\begin{equation}
\Vertpair{F} \leq \parenpair*{a^*}^{1/2}, \quad \Vertpair{J} \leq \parenpair*{a^*}^{1/2} \label{eq:fjsepnormbd}
\end{equation}

Moreover, the product $FJ$ is a Hilbert-Schmidt operator, or of the Schatten class $\schat{2}$ in a Hilbert space $\ell^2(\NN)$; see more in \cite{GohKri}, \cite{Simon}.  Indeed, 
\begin{align}
\anglepair*{FJe_k, e_m} & = a_k^{1/2} a_m^{1/2} f_m(y_k) \label{eq:fjprod}
\intertext{and}
\sum_{k, m = 1}^{\infty} \vertpair*{\anglepair*{FJ e_k, e_m}}^2 & = \sum_{k, m = 1}^{\infty} a_k a_m \vert f_m(y_k) \vert^2 \leq \parenpair*{a^*}^2
\end{align}
so $\snorm{FJ}{2} \leq a^*$.  

By H\"{o}lder inequality for Schatten classes (\cite{GohKri} or \cite{Simon}),
\begin{equation}
\snorm{BCD}{2/3} \leq \snorm{B}{2} \snorm{C}{2} \snorm{D}{2}
\end{equation}
so $(FJ)^3 \in \schat{2/3}$ and has a representation
\begin{equation}
(FJ)^3 = \sum_{k = 1}^{\infty} c_k \anglepair{\cdot, f_k} h_k, \quad c > 0, 
\end{equation}
where $\Vertpair{f_k}, \Vertpair{h_k} \leq 1$ and 
\begin{equation}
\sum_{k = 1}^{\infty} c_k^{2/3} < \infty.
\end{equation}
Therefore, 
\begin{equation}
A^4 = J(FJ)^3 F = \sum_{k = 1}^{\infty} c_k \anglepair{F(\cdot), f_k} Jh_k
\end{equation}
has $\disp \frac{2}{3}$-property \eqref{eq:twothirds} and 
\begin{equation} \label{eq:sumeigenq}
\sum_{j = 1}^{\infty} \vertpair*{\lambda_j(A^q)} < \infty \text{ for all } q \geq 4, 
\end{equation}
with 
\begin{equation} \label{eq:traq}
\tra A^q = \sum_{j = 1}^{\infty} \lambda_j(A^q)
\end{equation}
More careful geometric analysis, based on approximative characteristics of operators \cite{MitPel}, \cite{Piet} --- if we use \cite{Koniga}, or \cite[Theorem 4.a.6, p. 227]{Konigb}  --- shows that we can lower $q$ in \eqref{eq:sumeigenq}, \eqref{eq:traq} to $3$.  Indeed, $(FJ)^2$ is in $\schat{1}(\ell^2(\NN))$, so there are finite-dimensional operators $G_n$, $\ranky G_n \leq n$, such that 
\begin{equation}
\sum_n \alpha_n < \infty, \quad \text{ where }\alpha_n \declare \Vertpair{(FJ)^2 - G_n}
\end{equation}
Then 
\begin{equation}
\Vertpair{A^3 - J G_n F} \leq a^* \cdot \alpha_n
\end{equation}
and by \cite[Theorem 4.a.6]{Konigb} 
\begin{align}
A^3 \text{ is nuclear}, \\
\sum_j \vertpair{\lambda_j(A^3)} \leq 2 a^* \sum_n \alpha_n < \infty, \\
\intertext{and}
\tra A^3 = \sum_{j} \lambda_j(A^3).
\end{align}
But this remark will not improve our Theorem~\ref{thm:mainthm} (below) in an essential way (just in \eqref{eq:trabigoffres} we can say $p \geq p_* \geq 3q*$).  

In a Hilbert space $X = H$ by Lisdki\u{\i} Theorem \cite{Lidskii}, for any trace-class operator $C \in \schat{1}$, 
\begin{gather}
\sum_{j = 1}^{\infty} \vertpair*{\lambda_j(C)} < \infty \\
\intertext{and}
\tra C = \sum_{j = 1}^{\infty} \lambda_j(C).
\end{gather}
Maybe, talking just about nuclear operators, M. Zelikin considered in \cite[Thm. 2]{Zelikin} only Hilbert spaces.

Before stating our main result let us recall \cite[Chapter VII, Sections 3 and 4]{DunSchI} elements of Riesz theory of compact operators.  

If $T: X \to X$ is compact its spectrum $\spectr{T}$ is discrete with $0$ being the only accumulation point, and it has the following properties
\begin{enumerate}[label = (\roman*)]
\item for any $\rho > 0$,  
$\spectr{T} \cap \bracepair{z: \vert z \vert \geq \rho}$ is a finite set; 
\item if
\begin{equation}
\locrad{\alpha} = \frac{1}{2} \min \bracepair{\vertpair{\alpha - \lambda} : \lambda \in \spectr{T}, \lambda \neq \alpha}
\end{equation}
[so $\locrad{\alpha} > 0$ for any $\alpha \in \CC \setminus{0}$] and 
\begin{equation}
\locproj{\alpha} = \frac{1}{2\pi i} \int\limits_{\vertpair{z - \alpha} = \locrad{\alpha}} (z - T)^{-1} \, dz,
\end{equation}
then
\begin{equation}
\mcy{\alpha} = \ranky \locproj{\alpha} < \infty, \quad \alpha \in \CC \setminus \bracepair{0}
\end{equation}
with 
\begin{equation}
\mcy{\alpha} = 0 \quad \text{ if and only if } \alpha \not\in \spectr{T}.
\end{equation}
For $\alpha \in \spectr{T} \setminus{0}$, $m(\alpha)$ is an algebraic multiplicity of an eigenvalue $\alpha$.  
\end{enumerate}
The operational calculus \cite[Chapter VII, Sections 3 and 4]{DunSchI} explains that for any $\rho > 0$ such that 
\begin{equation}
\spectr{T} \cap \bracepair{ \vertpair{z} = \rho} = \emptyset
\end{equation}
we have 
\begin{equation}
T = \sum_{\vert \alpha \vert > \rho} T(\alpha) + S, \text{ where } T(\alpha)  = \frac{1}{2\pi i} \int\limits_{\vertpair{z - \alpha} = \locrad{\alpha}} z(z - T)^{-1} \, dz
\end{equation}
is an operator of rank $\mcy{\alpha}$ with
\begin{equation}
\spectr{T(\alpha)} = \bracepair{\alpha},
\end{equation}
and
\begin{equation}
S = \frac{1}{2\pi i} \int\limits_{\vert z \vert = \rho} z(z - T)^{-1} \, dz.
\end{equation}
Moreover, for any entire function $F(z)$, say, for polynomials, 
\begin{equation} \label{eq:fcnloc}
F(T) = \sum_{\vert \alpha \vert > \rho} F(T(\alpha)) + F(S),
\end{equation}
where by the Riesz-Cauchy formulae,
\begin{equation}
F(T(\alpha)) =  \frac{1}{2\pi i} \int\limits_{\vertpair{z - \alpha} = \locrad{\alpha}} F(z)(z - T)^{-1} \, dz, \quad F(S) = \frac{1}{2\pi i} \int\limits_{\vert z \vert = \rho} F(z)(z - T)^{-1} \, dz.
\end{equation}
It follows that
\begin{gather}
\tra F(T(\alpha))  = F(\alpha) \cdot \mcy{\alpha}. \\
F(T(\alpha))  = 0 \quad \text{ if } F^{(j)}(\alpha) = 0, \quad 0 \leq j \leq \mcy{\alpha}. \label{eq:zerocond}
\end{gather}
Now we are ready to prove
\begin{mainthm} \label{thm:mainthm}
Let $T$ be a compact operator in a Banach space $X$, and some power $T^{q_*}$ is a nuclear operator.  Then $\spectr{T}$ is $\zd{d}$--symmetric, i.e., for any $\beta \in \CC \setminus \bracepair{0}$, 
\begin{equation} \label{eq:symdef}
\begin{split}
\mcy{\beta \omega^k} = \mcy{\beta} \text{ for all } k = 0, 1, \dotsc d - 1, \quad \omega = \exp \parenpair*{i \frac{2\pi}{d}}
\end{split}
\end{equation}
if and only if 
\begin{equation} \label{eq:trabigoffres}
\tra T^{dp + r} = 0, \quad 1 \leq r \leq d -1 ,
\end{equation}
for all sufficiently large $p$, say $p \geq p_* \geq 4 q_*$.
\end{mainthm}
Of course, if $d = 2$, this is an extension of \cite{Zelikin} , Thm. 2, to a Banach case.
\begin{proof}[Proof.  Part 1: \eqref{eq:symdef} $\Rightarrow$ \eqref{eq:trabigoffres}]  This is an ``algebraic'' claim although first we notice: the assumption $p \geq 4q_*$ guarantees that all operators $T^n$, $n = dp + r$, in \eqref{eq:trabigoffres} satisfy $\frac{2}{3}$--condition so by Grothendieck theorem 
\begin{equation}
\tra T^n = \sum_{j = 1}^{\infty} \lambda_j(T^n)
\end{equation}
and the absolute convergence permits to rearrange the terms of the right sum as we wish to write
\begin{equation} \label{eq:tratnbetterwrite}
\tra T^n = \sum \mu \cdot \mcy{\mu; T^n}
\end{equation}
With
\begin{equation}
\mcy{\mu; T^n} = 0 \text{ for } \mu \not\in \spectr{T^n} 
\end{equation}
we can ``add'' the terms with $\mu \not\in \spectr{T^n}$ and this does not change the right side in \eqref{eq:tratnbetterwrite}.  For
\begin{equation} \label{eq:gcdform}
n = dp + r \in \eqref{eq:trabigoffres} \text{ define } g = \gcd \bracepair{r, d}
\end{equation}
so
\begin{equation} \label{eq:gcdrules}
r = ag, \quad d = bg, \quad (a, b) = 1
\end{equation}
and with $r \leq d - 1$ we have $1 \leq a < b$.
For any $\mu \in \CC \setminus \bracepair{0}$ take its $\zd{b}$-orbit, i.e.,
\begin{equation} \label{eq:orbitdef}
\widetilde{\mu} = \bracepair{\mu \cdot \tau^j: 0 \leq j < b}, \quad \tau = \omega^q = \exp \parenpair*{i \frac{2\pi}{b}}.
\end{equation}
The sum in \eqref{eq:tratnbetterwrite} could be written as
\begin{equation} \label{eq:tratnmajorrewrite}
\sum\limits_{\zd{b}-\text{orbits}} \, \, \sum_{j = 0}^{b-1} \mu \tau^j \cdot \mcy{\mu \tau^j; T^n}
\end{equation}
where for certainty $\mu$ in the orbit \eqref{eq:orbitdef} is chosen as $\mu = \vert \mu \vert e^{i \vartheta}$, $0 \leq \vartheta < \frac{2\pi}{b}$.  Now we will show that the sum in \eqref{eq:tratnmajorrewrite} over each orbit is equal to zero.  With numbers as in \eqref{eq:gcdform} put $\disp \kappa = \exp \parenpair*{i \frac{2\pi}{n}}$ so $\kappa^n = 1$ and notice that if $\mu = \lambda^n$, we choose
\begin{equation} \label{eq:npowerrewriteone}
\lambda = \vert \mu \vert^{1/n} e^{i \vartheta^{\prime}}, \quad \vartheta^{\prime} = \frac{\vartheta}{n},
\end{equation}
then 
\begin{equation} \label{eq:lampluspows}
\parenpair*{\lambda \omega^k}^n = \mu \omega^{k(dp + r)} = \mu \omega^{kr} = \mu \tau^{ak}
\end{equation}
and
\begin{equation} \label{eq:maineq}
\begin{split}
\sum_{j = 0}^{b - 1} \mu \tau^j \mcy{\mu \tau^j; T^n} &= \frac{1}{g} \sum_{k = 0}^{d-1} \mu \tau^{ak} \mcy{\mu \tau^{ak}; T^n} \eqprintnew{sign:first} \\
 & = \frac{1}{g} \sum_{k = 0}^{d-1} \left( \lambda \omega^k \right)^n \sum_{s = 0}^{n-1} \mcy{\lambda \omega^k \kappa^s; T} \eqprintnew{sign:second} \\
 & = \frac{1}{g} \sum_{s = 0}^{n-1} \sum_{k = 0}^{d-1} \left( \lambda \omega^k \right)^n  \mcy{\lambda \kappa^s \cdot \omega^k; T}  \eqprintnew{sign:third} \\
 & = \frac{1}{g} \sum_{s = 0}^{n-1} \mcy{\lambda \kappa^s; T} \mu \sum_{k = 0}^{d-1} \tau^{\alpha k} \eqprintnew{sign:fourth} \\
 &  \mu \left( \sum_{s = 0}^{n-1} \mcy{\lambda \kappa^s; T} \right) \left( \sum_{j = 0}^{b - 1} \tau^j \right) \eqprintnew{sign:fifth} 0 
\end{split}
\end{equation}
The steps in \eqref{eq:maineq} are justified in the following way.  \ref{sign:first} comes from \eqref{eq:lampluspows}.  \ref{sign:second} is just the change of order of the double summation.  \ref{sign:third} uses in essential way the theorem's assumption \eqref{eq:symdef} on $\mcy{\beta \omega^k}$ being independent on $k$.  \ref{sign:fourth} is bases on the properties of the roots $\omega, \, \tau$, $\omega^d = 1$, $\tau = \omega^g$ under \eqref{eq:gcdrules}.  Of course, in \ref{sign:fifth} $\disp \sum_{j = 0}^{b-1} \tau^j = 0$, and $\bracepair{\tau^{ak}}_{k = 0}^{d-1}$ runs $g$ times over $\bracepair{\tau^j}_{j = 0}^{b-1}$.  Part \eqref{eq:symdef} $\Rightarrow$ \eqref{eq:trabigoffres} is proven.
\end{proof}
\begin{proof}[Proof.  Part 2: \eqref{eq:trabigoffres} $\Rightarrow$ \eqref{eq:symdef}]  Take $\lambda \neq 0$ and as before 
\begin{equation}
n = dp_* + dp + r, \quad 1 \leq r \leq d - 1, \quad p \geq 0
\end{equation}
and $0 < \rho < \vertpair{\lambda}$ is such that
\begin{equation} \label{eq:separation}
\spectr{T} \cap \bracepair{z \in \CC: \vert z \vert = \rho} = \emptyset,
\end{equation}
with
\begin{equation}
\widetilde{\lambda} = \bracepair{\lambda \omega^k: 0 \leq k \leq d - 1}
\end{equation}
being the $\zd{d}$-orbit of $\lambda$.  Now we use \eqref{eq:fcnloc} for the special choice $F = F_{pr}$ with  
\begin{equation}
F_{pr}(z) = \parenpair*{\frac{z}{\lambda}}^{dp_* + dp + r} \varphi(z),
\end{equation}
where
\begin{align}
\varphi(z) &= \prod_{\substack{\vert \alpha \vert \geq \rho \\ \alpha \in \spectr{T} \\ \alpha \not\in \widetilde{\lambda}}} \parenpair*{\frac{z^d - \alpha^d}{\lambda^d - \alpha^d}}^{\mcy{\alpha}} =\\
& =  \psi(z^d), \quad \text{ and }\psi \text{ is a polynomial.}
\end{align}
Then by \eqref{eq:zerocond}
\begin{align}
\varphi(T(\alpha)) & = 0,\\
F_{pr}(T(\alpha)) & = 0, \quad \forall \alpha \not\in \widetilde{\lambda}, \quad \vert \alpha \vert > \rho
\end{align}
but for $\beta \in \widetilde{\lambda}$, i.e., $\beta = \lambda \omega^k$, 
\begin{equation}
\tra F_{pr}(T(\beta)) = \mcy{\beta} F_{pr}(\beta) = \mcy{\lambda \omega^k} \omega^{kr}.
\end{equation}
Therefore,
\begin{equation} \label{eq:fprtdecomprev}
\tra F_{pr}(T) = \sum_{k = 0}^{d-1} \omega^{kr} \mcy{\lambda \omega^k} + \tra F_{pr}(S)
\end{equation}
where
\begin{equation}
F_{pr}(S) = \left( \frac{T}{\lambda} \right)^{dp_*} \cdot \frac{1}{2\pi i} \int\limits_{\vert z \vert = \rho} \left( \frac{z}{\lambda} \right)^{dp + r} \varphi(z) (z - T)^{-1} \, dz.
\end{equation}
Put
\begin{equation}
\Phi = \max \bracepair{\vertpair{\varphi(z)}: \vertpair{z} \leq \rho}
\end{equation}
and with \eqref{eq:separation}
\begin{equation}
M = \max \bracepair{\Vertpair{R(z; T)}: \vert z \vert = \rho} < \infty
\end{equation}
Then 
\begin{equation} \label{eq:fprsnorm}
\snorm{F_{pr}(S)}{1} \leq C t^p, \quad \text{ any } r, \quad 1 \leq r \leq d - 1,
\end{equation}
where
\begin{equation} \label{eq:constdef}
C = \frac{\Phi \cdot M \cdot \rho \cdot \snorm{T^{dp_*}}{1}}{\vertpair{\lambda}^{dp*}} 
\end{equation}
and
\begin{equation} \label{eq:tbd}
t = \parenpair*{\frac{\rho}{\vert \lambda \vert}}^d < 1.
\end{equation}
Now by \eqref{eq:trabigoffres} and \eqref{eq:fprtdecomprev}
\begin{equation} \label{eq:zeroeq}
0 = \sum_{k = 0}^{d - 1} \omega^{kr} \mcy{\lambda \omega^k} + \xi_{pr} \quad \text{ for any } p \geq 1 \text{ and }r, \quad 1 \leq r \leq d - 1.
\end{equation}
The sum $\disp \sum_{k = 0}^{d-1}$ does not depend on $p$ but the remainder by \labelcref{eq:fprsnorm,eq:constdef,eq:tbd} have estimates
\begin{equation}
\vert \xi_{pr}\vert  \leq C t^p \quad \text{ so } \quad \xi_{pr} \to 0 \, \, (p \to \infty)
\end{equation}
This implies by \eqref{eq:zeroeq} 
\begin{equation} \label{eq:zerores}
\sum_{k = 0}^{d-1} \omega^{kr} \mcy{\lambda \omega^k} = 0, \quad \forall r, \quad 1 \leq r \leq d - 1
\end{equation}
or
\begin{equation}
y_k = \mcy{\lambda \omega^k}, \quad 1 \leq k \leq d - 1.
\end{equation}
is a solution of the system
\begin{equation} \label{eq:sysform}
\sum_{k = 1}^{d-1} \omega^{kr} y_k = - y_0, \quad 1 \leq r \leq d-1.
\end{equation}
Its determinant is of Vandermonde type so
\begin{equation}
\det \bracepair{\omega^{kr}}_{k, r = 1}^{d-1} \neq 0,
\end{equation}
and the identities
\begin{equation}
\sum_{k = 0}^{d-1} (\omega^r)^k = 0, \quad \forall r, \quad 1 \leq r \leq d - 1
\end{equation}
show that by \eqref{eq:sysform}
\begin{equation}
y_k = y_0, \text{ i.e., } \mcy{\lambda \omega^k} = \mcy{\omega}, \, \forall k, \quad 1 \leq k \leq d - 1.
\end{equation}
This proves that the multiplicity function $m$ is constant on $\zd{d}$--orbits in $\CC \setminus \bracepair{0}$, and \eqref{eq:symdef} is proven.
\end{proof}

It is worth to notice that the proof of Part II does not use any form of Grothendick or Lidskii thoerem but it uses only properties of a linear function \textit{Trace} on $\mathcal{N}(X; X)$ and an elementary formula for $\tra K$ when $K$ is an operator of finite rank.

\section*{Acknowledgements}  I thank M. Zelikin who told me his results \cite{Zelikin} in September 2007 during the 18th Crimean Autumn Mathematical School--Symposium, in Laspi--Batiliman, Crimea.  At that time I realized that his result could be extended to the Banach space case.  In the late 2007 I had useful discussions with P. Djakov, E. Gorin, H. K\"{o}nig, M. Solomyak.  Recently, during the Aleksander Pelczynski Memorial Conference, July 2014, Bedlewo, Poland, O. Reinov brought my attention to his papers \cite{Rei}, \cite{ReiLat}.  I thank C. Baker and O. Reinov for recent discussions.

\nocite{*}
\printbibliography

\end{document}